\documentclass[conference]{IEEEtran}
\IEEEoverridecommandlockouts
\usepackage{cite}
\usepackage{amsmath,amssymb,amsfonts}
\usepackage{algorithmic}
\usepackage{graphicx}
\usepackage{textcomp}
\usepackage{xcolor}
\usepackage{amsmath}
\usepackage{mathtools, stmaryrd}
\usepackage{amssymb}
\usepackage{graphicx}
\usepackage{subcaption}
\usepackage{url}
\usepackage{tabularx}
\usepackage{supertabular}
\usepackage{multirow}
\usepackage{boxedminipage}
\usepackage{eurosym}
\usepackage{enumitem}   
\usepackage[us]{datetime} 
\usepackage[bookmarks=false]{hyperref}
\hypersetup{
     colorlinks=true,
     linkcolor=blue,
     filecolor=blue,
     citecolor=black,      
     urlcolor=cyan,
     }
\usepackage{flushend}

\newcommand{\CREspecifications}{\formatdate{12}{07}{2019}}

\usepackage{fancyhdr}
\def\BibTeX{{\rm B\kern-.05em{\sc i\kern-.025em b}\kern-.08em
    T\kern-.1667em\lower.7ex\hbox{E}\kern-.125emX}}
\begin{document}

\title{Stochastic and deterministic formulations for capacity firming nominations\\
\thanks{The first two authors would like to thank John Cockerill and Nethys for their financial support.}
}

\author{\IEEEauthorblockN{Jonathan Dumas\IEEEauthorrefmark{1}, Bertrand~Corn\'elusse\IEEEauthorrefmark{1}, Antonello Giannitrapani\IEEEauthorrefmark{2}, Simone Paoletti\IEEEauthorrefmark{2}, Antonio Vicino\IEEEauthorrefmark{2}}
\IEEEauthorblockA{\IEEEauthorrefmark{1} \textit{Department of Electrical Engineering and Computer Science} \\
\textit{University of Li\`ege, Belgium}\\
\{jdumas, bertrand.cornelusse\}@uliege.be
}
\IEEEauthorblockA{\IEEEauthorrefmark{2} \textit{Dipartimento di Ingegneria dell'Informazione e Scienze Matematiche} \\
	\textit{Universit\`a di Siena, Italy}\\
	\{giannitrapani, paoletti, vicino\}@dii.unisi.it
}
}

\maketitle

\thispagestyle{fancyplain}

\begin{abstract}
This paper addresses the \textit{energy management} of a grid-connected photovoltaic plant coupled with a battery energy storage device, within the \textit{capacity firming} specifications of the French Energy Regulatory Commission. The paper contributions are positioned in the continuity of the studies adopting stochastic models for optimizing the bids of renewable energy sources in a day-ahead market by considering a storage device. The proposed deterministic and stochastic approaches are optimization problems formulated as quadratic problems with linear constraints. The case study is a real microgrid with PV production monitored on-site. The results demonstrate the validity of the stochastic formulation by using an ideal predictor that produces unbiased PV scenarios.
\end{abstract}


\newcommand{\SOC}{\ensuremath{s}}
\newcommand{\charge}{\ensuremath{P^\text{cha}}}
\newcommand{\discharge}{\ensuremath{P^\text{dis}}}
\newcommand{\maxcharge}{\ensuremath{\overline{S}}}
\newcommand{\mincharge}{\ensuremath{\underline{S}}}
\newcommand{\chargerate}{\ensuremath{\overline{S_P}}}
\newcommand{\dischargerate}{\ensuremath{\underline{S_P}}}
\newcommand{\retentionRate}{\ensuremath{\eta^{\text{retention}}}}
\newcommand{\chargeEfficiency}{\ensuremath{\eta^{\text{cha}}}}
\newcommand{\dischargeEfficiency}{\ensuremath{\eta^{\text{dis}}}}
\newcommand{\initialCharge}{\ensuremath{S^\text{init}}}
\newcommand{\minEndCharge}{\ensuremath{\underline{S}^\text{end}}}
\newcommand{\maxEndCharge}{\ensuremath{\overline{S}^\text{end}}}
\newcommand{\finalCharge}{\ensuremath{S^{\text{end}}}}

\newcommand{\maxExportToGrid}{\ensuremath{E}_{t}^{\text{cap}}}

\newcommand{\price}[1]{\ensuremath{\pi^{\text{#1}}}}
\newcommand{\gridSalePrice}{\price{exp}}

\newcommand{\exportrealizations}{\ensuremath{e}^\text{m}}

\newcommand{\nominations}{\ensuremath{e}^{\star}}
\newcommand{\nominationsvector}{\ensuremath{\underline{e}}^{\star}}

\newcommand{\export}{\ensuremath{e}}
\newcommand{\exportInit}{\ensuremath{e}^\text{ini}}
\newcommand{\PositiveEnergyDeviation}{\ensuremath{\Delta e^+}}
\newcommand{\NegativeEnergyDeviation}{\ensuremath{\Delta e^-}}
\newcommand{\nonSteerable}{\ensuremath{P^\text{m}}}
\newcommand{\nonSteerableMax}{\ensuremath{\underline{P}_c}}
\newcommand{\production} {\ensuremath{P}}
\newcommand{\productionForecast} {\ensuremath{\widehat{P}}}

\newcommand{\temperature}{\ensuremath{T}}
\newcommand{\irradiance}{\ensuremath{I}}
\newcommand{\clearskyirradiance}{\ensuremath{I}^\text{cs}}
\newcommand{\clearskyirradianceForecast}{\ensuremath{\widehat{I}}^\text{cs}}
\newcommand{\irradianceForecast}{\ensuremath{\widehat{I}}}
\newcommand{\temperatureForecast}{\ensuremath{\widehat{T}}}

\newcommand{\foralltmarketperiod}{\ensuremath{\forall t \in \mathcal{T}}}
\newcommand{\forallw}{\ensuremath{\forall \omega \in \Omega}}
\newcommand{\foralltrealtimeperiod}{\ensuremath{\forall t' \in \mathcal{T'}}}

\section{Introduction}
The capacity firming framework is mainly designed for islands or isolated markets. For instance, the French Energy Regulatory Commission (CRE) publishes capacity firming tenders and specifications\footnote{\url{https://www.cre.fr/}.}. The system considered is a grid-connected photovoltaic (PV) plant with a battery energy storage system (BESS) for firming the PV generation. At the tendering stage, offers are selected on the electricity selling price\footnote{Other less important criteria may count in the tender decision.}. Then, the successful tenderer builds its plant and sells the electricity exported to the grid at the contracted selling price, but according to a well-defined daily nomination and penalization scheme. The electricity to be injected in or withdrawn from (for auxiliaries) the grid must be nominated the day-ahead, and nominations must satisfy ramping power constraints. Under specific conditions, it is possible to renominate within the day at specified market gates. The remuneration is calculated a \textit{posteriori} by multiplying the realized exports by the contracted selling price minus a penalty. The deviations of the realized exports from the nominations are penalized through a function defined in the specifications of the tender.

The optimal day-ahead bidding strategies of a plant composed of only a production device have been addressed in, \textit{e.g.},  \cite{pinson2007trading,bitar2012bringing,giannitrapani2014bidding,giannitrapani2015bidding}. The optimal offer turns out to be a suitable percentile of the PV/wind power cumulative distribution function. Under the assumption of time-invariant statistics of power generation, the cumulative distribution functions can be estimated from historical data of the power generated by the plant. This assumption is not always justified especially for PV power generation. In \cite{giannitrapani2015bidding}, the authors investigate two approaches to properly take into account the effects of seasonal variation and non-stationary nature of PV power generation in the estimation of PV power statistics. However, incorporating energy storage in the framework is still an open problem, and the literature provides several approaches and methodologies to this end.
An optimal power management mechanism for a grid-connected PV system with storage is implemented in \cite{riffonneau2011optimal} using Dynamic Programming (DP) and is compared with simple ruled-based management. The sizing and control of an energy storage system to mitigate wind power uncertainty is addressed by \cite{haessig2014dimensionnement,haessig2013aging,haessig2015energy} using stochastic dynamic programming (SDP). The framework is similar to the CRE PV capacity firming tender with a wind farm operator that is committed on a day-ahead basis to a production engagement. Finally, three distinct optimization strategies, mixed-integer quadratic programming, simulation-based genetic algorithm, and expert-based heuristic are empirically compared by \cite{n2019optimal} in the CRE framework. 

This paper addresses the energy management of a grid-connected PV plant and BESS. This topic is studied within the capacity firming specifications of the CRE, in line with the tender AO-CRE-ZNI 2019 published on $\CREspecifications$, using the MiRIS microgrid case study.
The capacity firming problem can be decomposed in two steps. The first step consists of computing the day-ahead nominations. The second step consists of computing the renominations and the set-points in real-time to minimize the energy and ramp power deviations from nominations. This paper focuses on the first step and proposes both a stochastic and a deterministic formulation. The main goal of this paper is to validate the stochastic approach by using an ideal predictor providing unbiased PV scenarios. Thus, the BESS efficiencies are perfect and the degradation is not taken into account for the sake of simplicity. Different levels of prediction accuracy are evaluated. Then, the results are compared with those of the deterministic formulation, assuming perfect forecasts returned by an oracle. 
Both deterministic and stochastic approaches result in optimization problems formulated as quadratic problems with linear constraints. The considered case study is a real microgrid with PV production monitored on-site.

The paper is organized as follows. Section \ref{sec:notation} provides the notation. Section \ref{sec:nominations_process} details the day-ahead nomination process. Section \ref{sec:Problem formulation} proposes the deterministic and stochastic formulations of the nomination process. Section \ref{sec:case_study} introduces the MiRIS microgrid case study and presents the results.
Conclusions are drawn in Section \ref{sec:conclusions}. Annex \ref{sec:scenario_generation_methodology} describes the methodology to generate the set of unbiased PV scenarios.\\

\section{Notation}\label{sec:notation}

\subsection{Sets and indices}
\begin{supertabular}{l p{0.8\columnwidth}}
	Name & Description \\
	\hline
	$t$ & Market period index. \\
	$\mathcal{T}$ & $ \{1,2, \ldots, T\}$ set of market periods in a day. \\
	$\mathcal{D}$ & $ \cup_{i=1}^{i=D} \mathcal{T}_i$ set of $D$ days of market periods. \\
	$\omega$ & PV scenario index. \\
	$\Omega$ & $\{\omega_1,...,\omega_i, ..., \omega_{|\Omega|}\}$ set of PV scenarios. \\
\end{supertabular}

\subsection{Parameters}
\begin{supertabular}{l p{0.6\columnwidth} l}
	Name & Description & Unit \\
	\hline
	$\exportrealizations_t$ & Measured exports at $t$. & kWh\\
	$\nonSteerable_t$ & Measured power generation at $t$. & kW \\
	$\nonSteerableMax$ & Installed PV capacity. & kWp \\
	$\productionForecast_t$ & PV point forecast at $t$. & kW \\
	$\productionForecast_{t,\omega}$ & PV scenario $\omega$ at $t$. & kW \\
	$p_\omega$ & Probability of PV scenario $\omega$ &  /\\
	$\dischargerate$, $\chargerate$ & BESS maximum (dis)charging power. & kW \\
	$\mincharge$, $\maxcharge$ & Minimum/maximum capacity of the BESS.& kWh \\
	$\dischargeEfficiency$, $\chargeEfficiency$ & (Dis)Charging efficiency of the BESS. & /\\
	$\initialCharge$, $\finalCharge$ & Initial/final state of charge of the BESS. & kWh \\
	$\maxExportToGrid$   & Maximum export power to the grid. & kW \\
	$\gridSalePrice_t$ & Contracted selling price at $t$. & \euro/kWh \\
	$\pi^e$ & Slack price of $f^e$. & $\frac{\text{\euro}}{\text{kWh}^2}$ \\
	$\pi_{\maxcharge}$ & BESS CAPEX price. & \euro/kWh \\
	$\Delta_\tau$ & Market period duration. & minutes \\
	$\Delta E$ & Energy deadband between nomination and export. & kWh \\
	$\Delta P^{\star}$  & Power deadband between two consecutive nominations. & kW \\
	$f^e$ & Quadratic penalty function. & \euro \\
	
\end{supertabular}

\subsection{Variables}
\noindent The index $\omega$ is omitted in the deterministic case.
\begin{supertabular}{l l p{0.5\columnwidth} l}
	Name & Range & Description & Unit \\
	\hline
	$\production_{t,\omega} $ & $ [0,\productionForecast_t]$ & PV generated at $t$ in scenario $\omega$. & kW \\
	$\charge_{t,\omega} $ & $ [0,\chargerate]$ & Charging power used at $t$ in scenario $\omega$.  & kW \\
	$\discharge_{t,\omega} $ & $ [0, \dischargerate]$ & Discharging power used at $t$ in scenario $\omega$. & kW \\
	$\nominations_t$ & $ \mathbb{R}_+$ & Nomination at $t$. & kWh\\
	$\export_{t,\omega}$ & $ \mathbb{R}_+$ & Export  at $t$ in scenario $\omega$. & kWh\\
	$\PositiveEnergyDeviation_{t,\omega}$ & $ \mathbb{R}_+$ & Positive deviation between export and nomination at $t$ in scenario $\omega$.  & kWh\\
	$\NegativeEnergyDeviation_{t,\omega}$ & $ \mathbb{R}_+$ & Negative deviation between export and nomination  at $t$ in scenario $\omega$.  & kWh\\
	$\SOC_{t,\omega}$ & $ [\mincharge, \maxcharge]$ & State of charge of the BESS  at $t$ in scenario $\omega$. & kWh\\
	
\end{supertabular}

\section{Day-ahead nomination process}\label{sec:nominations_process}

\begin{figure}[tb]
	\centering
	\includegraphics[width=\linewidth]{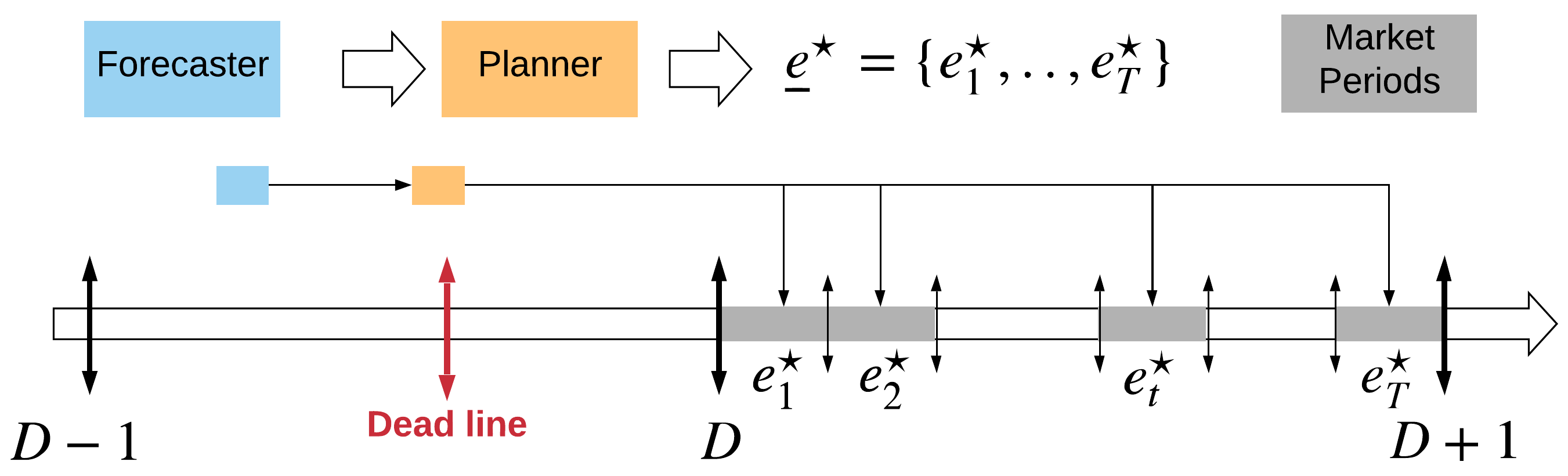}
	\caption{Day-ahead nomination process.}
	\label{fig:capacity_firming_nomination_process}
\end{figure}

Figure \ref{fig:capacity_firming_nomination_process} illustrates the day-ahead nomination process.
Each day is composed of $T$ market periods. The set of market periods is denoted by $\mathcal{T}$, and the market period duration is a constant value $\Delta_\tau$. 
The planner computes on a day-ahead basis a vector of nominations $\nominationsvector= [ \nominations_1,\cdots,\nominations_T ]^\intercal$, composed of $T$ values, based on forecasts of PV generation. The nominations are accepted if they satisfy the ramping power constraints
\begin{align}\label{eq:nominations_ramping_power_constraint}	
\frac{|\nominations_t-\nominations_{t-1}|}{\Delta_\tau} & \leq  \Delta P^{\star},  \ \ \foralltmarketperiod \setminus \{1\} ,
\end{align}
with $\Delta P^{\star}$ being a ramp limit, a fraction of the total installed capacity $\nonSteerableMax$, determined at the tendering stage, in our case of interest.
For a given market period, the net remuneration $r^n_t$ of the plant is proportional to the export $\exportrealizations_t$ minus a penalty $f^e(\nominations_t,\exportrealizations_t)$, with $\gridSalePrice_t$ being the contracted selling price
\begin{align}\label{eq:remuneration}	
r^n_t = \gridSalePrice_t  \exportrealizations_t  - f^e(\nominations_t,\exportrealizations_t),  \ \foralltmarketperiod.
\end{align}
The penalty function $f^e$ depends on the specifications of the tender. In this study, $f^e$ is approximated as  
\begin{align}\label{eq:energy_penalization}
f^e(\nominations_t,\exportrealizations_t)  & = \pi^e \bigg(  \max \big(0,|\nominations_t-\exportrealizations_t|- \Delta E \big) \bigg)^2
\end{align} 
where the deadband $\Delta E$ is a fixed fraction of the total installed capacity, and $\pi^e$ is a slack price (\euro/$\text{kWh}^2$).

\section{Problem formulation}\label{sec:Problem formulation}

Deterministic (D) and stochastic (S) formulations of the day-ahead nomination problem are compared. The deterministic formulation is used as a reference to validate the stochastic approach by considering perfect knowledge of the future (D$^\star$). Both approaches consider only exports to the grid\footnote{The imports from the grid are allowed only under specific conditions into the contract.}. The optimization variables and the parameters are defined in Section~\ref{sec:notation}.

\subsection{Deterministic formulation}\label{sec:deterministic_formulation}

The quadratic energy penalty $f^e$, defined in (\ref{eq:energy_penalization}), is applied to the deviations from the nominations. This formulation considers the nominations $\nominations_t$ and the exports $\export_t$ in the objective function $J_D$ to minimize
\begin{align}\label{eq:D_objective_1}	
J_D &  =  \sum_{t\in \mathcal{T}} - \gridSalePrice_t  \export_t +  f^e(\nominations_t, \export_t), 
\end{align}
which can be rewritten as
\begin{align}\label{eq:D_objective_2}	
J_D &  =\sum_{t\in \mathcal{T}} - \gridSalePrice_t  \export_t + \pi^e  \big( {(\PositiveEnergyDeviation_t)}^2 +  {(\NegativeEnergyDeviation_t)}^2 \big) 
\end{align}
by introducing the variables $\PositiveEnergyDeviation_t \in \mathbb{R}_+$ and $\NegativeEnergyDeviation_t \in \mathbb{R}_+$ defined as 
\begin{subequations}
	\label{eq:D_energy_deviation_constraints}	
	\begin{align}
	- \PositiveEnergyDeviation_t & \leq  - ( \export_t - \nominations_t- \Delta E  ), \ \foralltmarketperiod \\
	- \NegativeEnergyDeviation_t & \leq    - ( \nominations_t-\export_t -\Delta E ), \ \foralltmarketperiod .
	\end{align}
\end{subequations}
The optimization variables are $\nominations_t$, $\export_t$, $\PositiveEnergyDeviation_t$, $\NegativeEnergyDeviation_t$, $\production_t$, $\charge_t$, $\discharge_t$, and $\SOC_t$ (cf. Section~\ref{sec:notation}). 
From (\ref{eq:nominations_ramping_power_constraint}), the ramping power constraints are
\begin{subequations}
	\label{eq:D_nominations_ramping_constraints}	
	\begin{align}
	\frac{\nominations_t-\nominations_{t-1}}{\Delta_\tau} & \leq  \Delta P^{\star}, \ \foralltmarketperiod \setminus \{1\} \\
	\frac{\nominations_{t-1}-\nominations_t}{\Delta_\tau} & \leq   \Delta P^{\star}, \ \foralltmarketperiod \setminus \{1\}.
	\end{align}
\end{subequations}
The ramping constraint on $\nominations_1$ is deactivated to decouple consecutive days of simulation.
The set of constraints that bound the $\production_t$, $\charge_t$, $\discharge_t$, and $\SOC_t$ variables are, $\foralltmarketperiod$
\begin{subequations}
	\label{eq:D_action_constraint_set}	
	\begin{align}
	&\production_t \leq \productionForecast_t, \quad 
	\charge_t \leq \chargerate, \quad
	\discharge_t \leq \dischargerate    \\
	& \mincharge \leq \SOC_t \leq \maxcharge 
	\end{align}
\end{subequations}
where $\productionForecast_t$ are point forecasts of PV power generation. The power balance equation and the export constraints are
\begin{subequations}\label{eq:D_energy_flows}	
	\begin{align}
	\export_t/\Delta_\tau &  = \production_t +   (\discharge_t - \charge_t ), \ \foralltmarketperiod  \\
	\export_t/\Delta_\tau &  \leq \maxExportToGrid, \ \foralltmarketperiod \\
	\nominations_t/\Delta_\tau &  \leq \maxExportToGrid, \ \foralltmarketperiod.
	\end{align}
\end{subequations}
The dynamics of the BESS state of charge is
\begin{subequations}
	\label{eq:D_soc_dynamic}
	\begin{align}
	&\SOC_1 - \Delta_\tau (  \chargeEfficiency  \charge_1 - \frac{\discharge_1}{\dischargeEfficiency} ) = \initialCharge \\
	&\SOC_t - \SOC_{t-1} - \Delta_\tau (  \chargeEfficiency  \charge_t - \frac{\discharge_t}{\dischargeEfficiency}  )=0 \quad \foralltmarketperiod \setminus \{1\}\\
	&\SOC_T = \finalCharge = \initialCharge.
	\end{align}
\end{subequations}
The parameters $\finalCharge$ and $\initialCharge$ are introduced to decouple consecutive days of simulation.

\subsection{Deterministic formulation with perfect forecasts}\label{sec:deterministic_formulation_perfect}

With perfect forecasts, the above problem becomes
\begin{align}\label{eq:D_perfect_objective}	
\min & \ J_{D^\star}  =  \sum_{t\in \mathcal{T}} - \gridSalePrice_t  \export_t +  \pi^e  \big( {(\PositiveEnergyDeviation_t)}^2 +  {(\NegativeEnergyDeviation_t)}^2 \big) \\
& \ \text{s.t. } (\ref{eq:D_energy_deviation_constraints})-(\ref{eq:D_soc_dynamic}) \text{, with } \productionForecast_t = \nonSteerable_t \text{ in } (\ref{eq:D_action_constraint_set}).  \notag
\end{align}

\subsection{Stochastic formulation using a scenario-based approach}\label{sec:srochastic_formulation}

In the stochastic formulation of the day-ahead nomination problem, the objective is given by 
\begin{align}\label{eq:S_objective_1}	
J_S &  = \mathop{\mathbb{E}} \bigg[ \sum_{t\in \mathcal{T}} - \gridSalePrice_t \export_t +  f^e(\nominations_t, \export_t) \bigg] 
\end{align}
where the expectation is taken with respect to $\hat{P}_t$.
Using a scenario-based approach, (\ref{eq:S_objective_1}) is approximated by
\begin{align}\label{eq:S_objective_2}	
J_S &  =  \sum_{\omega \in \Omega} p_\omega  \sum_{t\in \mathcal{T}} \bigg[- \gridSalePrice_t \export_{t,\omega} +  f^e(\nominations_t, \export_{t,\omega})  \bigg]
\end{align}
with $p_\omega$ the probability of scenario $\omega\in \Omega$, and $\sum_{\omega \in \Omega} p_\omega = 1$. Then, by introducing $\PositiveEnergyDeviation_{t,\omega}$ and $\NegativeEnergyDeviation_{t,\omega}$ for each scenario $\omega$,  the problem becomes to solve
\begin{align}\label{eq:S_objective_3}	
  \min & \sum_{\omega \in \Omega} p_\omega  \sum_{t\in \mathcal{T}} \bigg[ - \gridSalePrice_t \export_{t,\omega}  +  \pi^e \big ( {\PositiveEnergyDeviation_{t,\omega}}^2 +  {\NegativeEnergyDeviation_{t,\omega}}^2 \big)   \bigg] \\
  & \text{ s.t } (\ref{eq:D_energy_deviation_constraints})-(\ref{eq:D_soc_dynamic}) \text{, with } \productionForecast_{t,\omega} \text{ instead of } \productionForecast_t \text{ in } (\ref{eq:D_action_constraint_set}), \ \forall \omega \in \Omega. \notag
\end{align}
All the optimization variables but $\nominations_t$ are now defined $\forall \omega \in \Omega$.

\subsection{Evaluation methodology}\label{sec:evaluation}

The second step of capacity firming, \textit{i.e.}, computing the set points in real-time, is required to assess the quality of the nomination process. 
However, since this paper focuses on the computation of day-ahead nominations, we simulate the second step with an ideal real-time controller\footnote{Using a real-time controller with intraday forecasts is required to assess the planner-controller. However, this study focus only on the nomination step.} once the nominations are fixed. The methodology to assess the nominations consists of minimizing 
\begin{align}\label{eq:eval_objective}	
J^\text{eval} &  =  \sum_{t\in \mathcal{T}}  - \gridSalePrice_t  \export_t  + f^e(\nominations_t, \export_t)
\end{align}
s.t (\ref{eq:D_energy_deviation_constraints})-(\ref{eq:D_soc_dynamic}) with $\productionForecast_t = \nonSteerable_t$ in (\ref{eq:D_action_constraint_set}) and given nominations $\nominations_t$ previously computed by the planner S. The optimization variables of (\ref{eq:eval_objective}) are $\export_t$, $\PositiveEnergyDeviation_t$, $\NegativeEnergyDeviation_t$, $\production_t$, $\charge_t$, $\discharge_t$, and $\SOC_t$. 
The optimal value of $J^\text{eval}_S$ is compared with the optimal value of $J_{D^\star}$ in (\ref{eq:D_perfect_objective}).

\section{MiRIS microgrid case study}\label{sec:case_study}

\begin{figure}[tb]
	\centering
	\includegraphics[width=\linewidth]{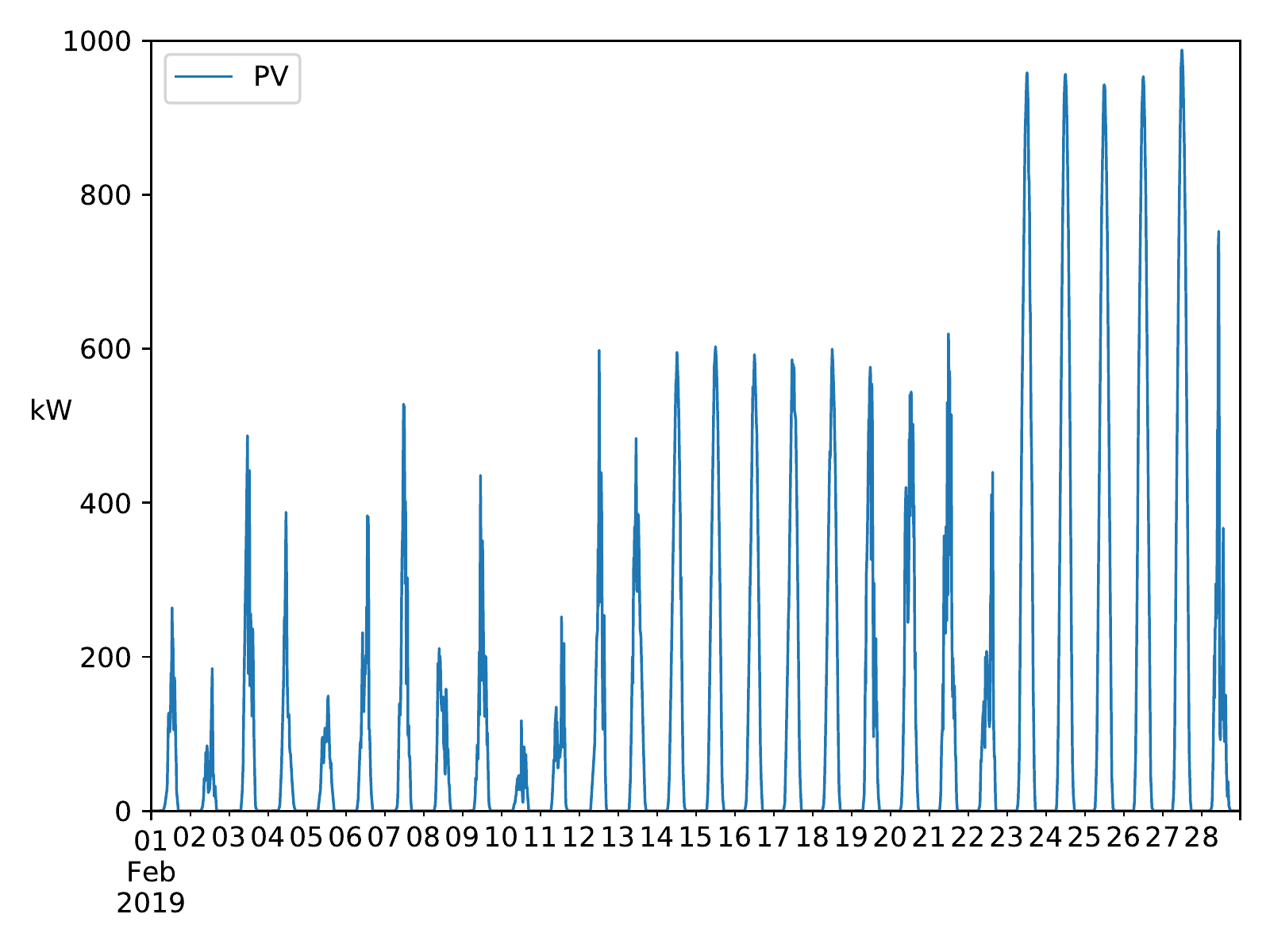}
	\caption{MiRIS February 2019 PV production.}
	\label{fig:pv_production}
\end{figure}

The MiRIS\footnote{\url{https://johncockerill.com/fr/energy/stockage-denergie/}} microgrid case study, located at the John Cockerill Group’s international headquarters in Seraing, Belgium, is composed of a PV production plant, a BESS, and a load. For the need of this study, only historical data of PV generation are required. The BESS capacity $\maxcharge$ is 1000 kWh, and the total PV installed capacity $\nonSteerableMax$ is 2000 kWp. The market period duration $\Delta_\tau$ is 15 minutes. The simulation dataset $\mathcal{D}$ is the month of February 2019. Figure~\ref{fig:pv_production} illustrates the MiRIS PV production and Table~\ref{tab:dataset_statistics} provides some key statistics.
Table~\ref{tab:indicators} defines the indicators used in this section.
\begin{table}[!htb]
	\caption{MiRIS February 2019 dataset statistics.}\label{tab:dataset_statistics}
	\centering
	\begin{tabular}{rrrrrr} \hline \hline
		$\nonSteerableMax$	& $ \overline{\nonSteerable}$ & $\nonSteerable_{\text{std}}$ & $\nonSteerable_{\max} $ & $\nonSteerable_{\%,\max} $ & $\nonSteerable_{\text{tot}}$ \\ \hline
		2000 & 104.6  & 202.4 & 988.1 & 49.4 & 70.3   \\ \hline \hline
	\end{tabular}
\end{table}
\begin{table}
	\caption{Comparison indicators.}\label{tab:indicators}
	\begin{tabularx}{\linewidth}{l X l}
		\hline \hline
		Name & Description & Unit \\
		\hline
		$\overline{\nonSteerable}$ & Averaged power. & kW\\
		$\nonSteerable_{\text{std}}$ & Power standard deviation. & kW \\
		$\nonSteerable_{\max}$ & Maximum power. & kW \\
		$\nonSteerable_{\text{tot}}$ & Total energy produced. & MWh \\
		$\nonSteerable_{\%,\max}$ & $\nonSteerable_{\max}$ divided by the total installed PV capacity $\nonSteerableMax$. & \% \\
		$\overline{t}_\text{CPU}$ & Averaged computation time per optimization problem. & s\\
		$[x]^D $ & Total of a variable $x_t$: $ \sum_{t\in \mathcal{D}} x_t$. & $x$ unit\\
		$[\production]^D$ & Total production. & MWh\\
		$\production_\%  $ & Production ratio: $\frac{[\production]^D }{[\nonSteerable]^D}$. & \% \\
		$\charge_\% $ & Charge ratio: $\frac{[\charge]^D }{[\production]^D}$. & \% \\
		$\maxcharge_\% $ & Percentage of days of the dataset where the BESS achieved its maximum storage capacity. & \% \\
		$\export_\% $ & Export ratio: $\frac{[\export]^D }{[\nominations]^D}$. & \% \\
		$R_{\max}$ & Maximum achievable revenue: $\gridSalePrice  [\nonSteerable]^D$. & k\euro \\
		$R^e$ & Gross revenue: $\gridSalePrice [\export]^D$. & k\euro \\
		$r^e $ & Maximum achievable revenue ratio: $ \frac{R^e}{R_{\max}}$. & \% \\
		$C^e $ & Quadratic penalty: $[f^e]^D$. & k\euro \\
		$R^{n,e}$ & Net revenue with quadratic penalty: $R^e - C^e$. & k\euro \\
		\hline \hline
	\end{tabularx}
\end{table}

It is of paramount importance to notice that the results of this case study are only valid for this dataset and cannot be extrapolated over an entire year without caution.
CPLEX\footnote{\url{https://www.ibm.com/products/ilog-cplex-optimization-studio}} 12.9 is used to solve all the optimization problems, on an Intel Core i7-8700 3.20 GHz based computer with 12 threads and 32 GB of RAM. Tables~\ref{tab:case_study_parameters} and~\ref{tab:battery_parameters} provide the case study and BESS parameters.
\begin{table}[tb]
	\caption{Case study parameters.}\label{tab:case_study_parameters}
	\centering
	\begin{tabular}{rrrrrr} \hline  \hline
		$\gridSalePrice$ & $\pi^e$& $\Delta_\tau$	& $\Delta P^{\star}$& $\maxExportToGrid$ & $ \Delta E$\\ \hline
		0.045  & 0.0045  & 15 & 10  & 2000 & 25 \\ \hline   \hline
	\end{tabular}
\end{table}
\begin{table}[tb]
	\caption{BESS parameters.}\label{tab:battery_parameters}
	\centering
	\begin{tabular}{rrrrrrrr} \hline  \hline
		$\maxcharge$ 	& $\mincharge$  & $\chargerate$&  $\dischargerate$  & $\chargeEfficiency$ & $\dischargeEfficiency$ & $\initialCharge$ & $\finalCharge$  \\ \hline
		1000 & 0 & 1000 & 1000 & 1 & 1 & 0 & 0 \\ \hline  \hline
	\end{tabular}
\end{table}

\subsection{Results for unbiased PV scenarios with fixed variance}

\begin{figure}[tb]
	\centering
	\begin{subfigure}{.25\textwidth}
		\centering
		\includegraphics[width=\linewidth]{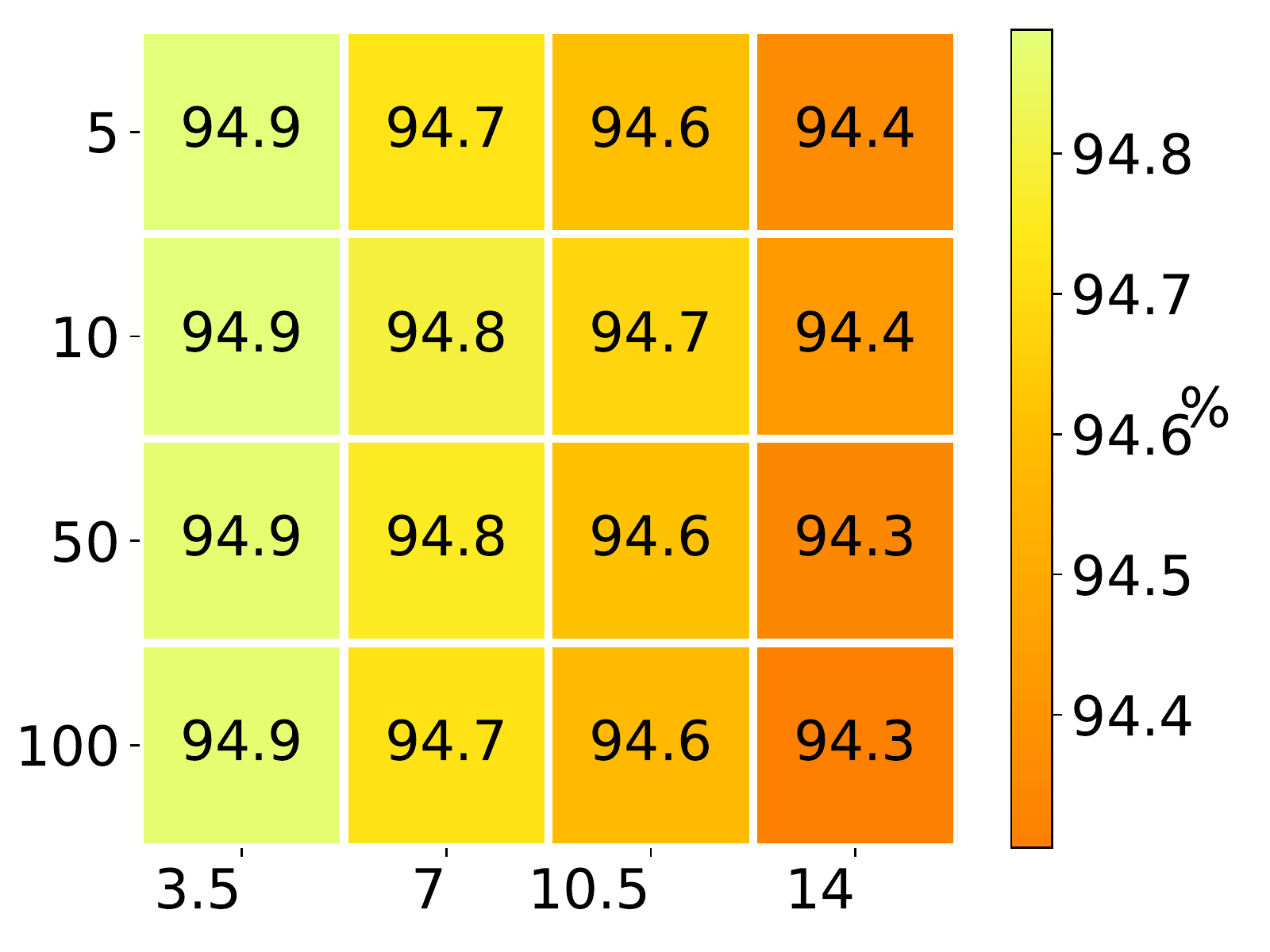}
		\caption{$\production_\%$.}
	\end{subfigure}%
	\begin{subfigure}{.25\textwidth}
		\centering
		\includegraphics[width=\linewidth]{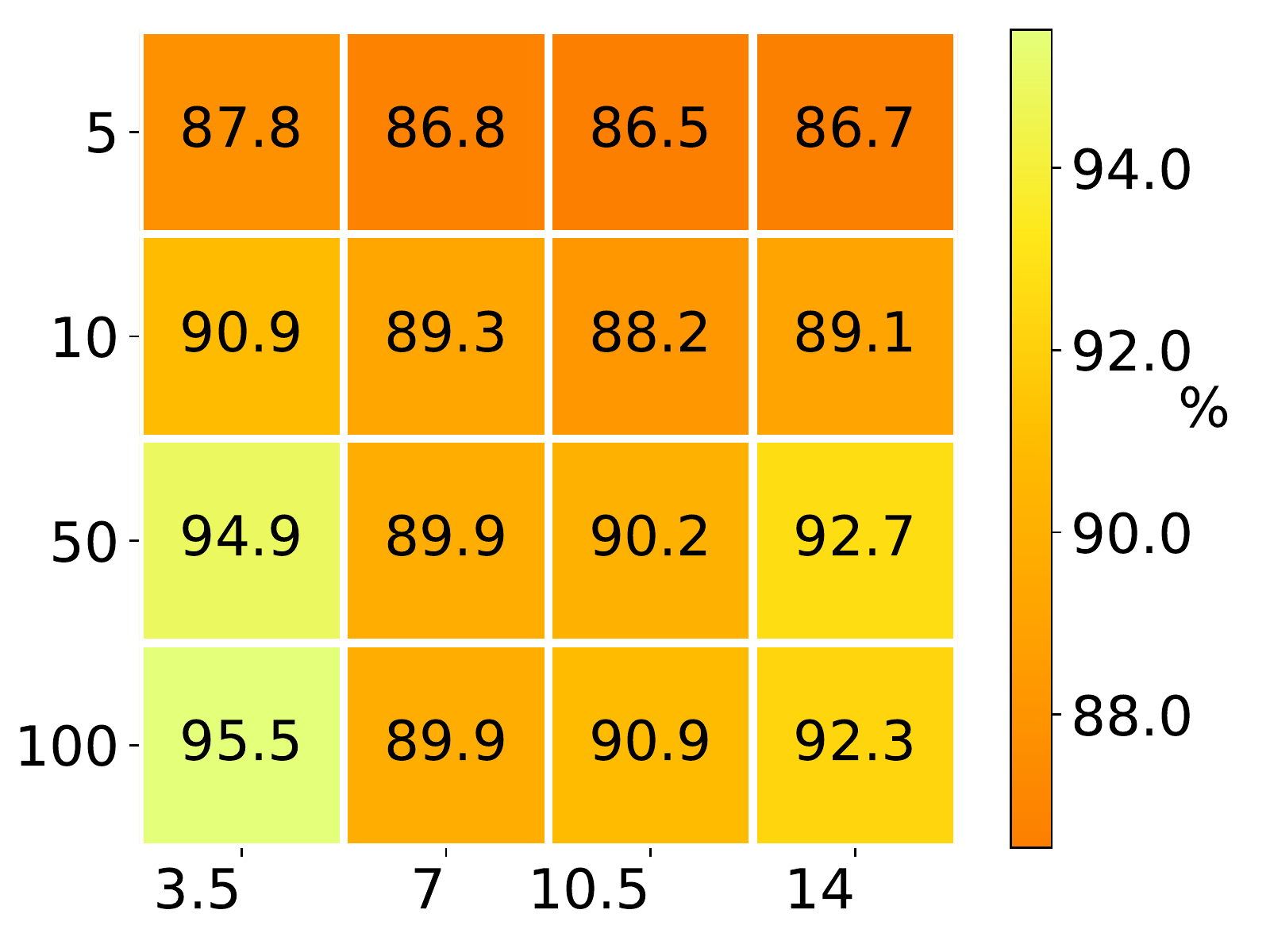}
		\caption{$\export_\%$.}
	\end{subfigure}
	\caption{Planner S ratio indicator with $|\Omega|=5, 10, 50, 100$, 	$\sigma = 3.5, 7, 10.5, 14\%$.}
	\label{fig:ratio_indicators_S}
\end{figure}

A set of unbiased PV scenarios is generated for several values of the standard deviation $\sigma$ of the prediction error. Table~\ref{tab:scenario_generation_parameters} shows the considered values of $\sigma$, expressed as a fraction of the actual PV generation. Moreover, Table~\ref{tab:scenario_generation_parameters} reports the cardinality of the generated scenario sets.
Table~\ref{tab:time_computation} compares the average computation time per optimization problem between planners S and D$^\star$. Note, The optimization problem of planner S with $|\Omega| = 1$ has the same number of variables and constraints as the  planner D$^\star$. The computation time is compatible with a day-ahead process even with 100 scenarios as it takes on average 7 seconds to compute the nominations for the day-ahead. 
Table~\ref{tab:ratio_indicators_D} and Figure~\ref{fig:ratio_indicators_S} provide the results of the ratio indicators, respectively, for the planners D$^\star$ and S.
\begin{table}[tb]
	\caption{Scenario generation parameters.}\label{tab:scenario_generation_parameters}
	\centering
	\begin{tabular}{rrrrr} \hline  \hline
		$\sigma$ & 3.5\% & 7\% & 10.5\% & 14\%\\ \hline
		$|\Omega|$ & 5& 10 & 50 & 100 \\ \hline  \hline
	\end{tabular}
\end{table}
\begin{table}[tb]
	\caption{Averaged computation times.}\label{tab:time_computation}
	\centering
	\begin{tabular}{lrrrrr} \hline  \hline
		$|\Omega|$ & 1 & 5& 10 & 50 & 100 \\ \hline
		\# variables& 769& 3457 & 6817  & 33697& 67297\\ \hline
		\# constraints& 1248& 5092 & 9897  & 48337& 96387\\ \hline
		$\overline{t}_\text{CPU}$ S& - & 0.3 & 0.8  & 3& 7\\ \hline
		$\overline{t}_\text{CPU}$ D$^\star$&  0.1& - & -  & -& -\\ \hline \hline
	\end{tabular}
\end{table}
\begin{table}[tb]
	\caption{Planner D$^\star$ ratio indicators.}
	\label{tab:ratio_indicators_D}
	\centering
	\begin{tabular}{rrrrr} \hline \hline
		$[\production]^D$ & $\production_\% $& $\charge_\%  $ & $\maxcharge_\%$ & $\export_\% $	\\ \hline
		66.7 & 94.9 &  29.6  & 17.9 & 76.2 \\ \hline \hline
	\end{tabular}
\end{table}

For all indicators, the results of both planners are almost equal with the smaller value of $\sigma$ and the highest value of $|\Omega|$, as expected. On average the curtailment of PV generation equals $5 \%$. The maximum $\export_\% $ value is achieved with $\sigma = 3.5\%$ because the nominations are more conservative when the variance increases, leading to a smaller ratio. On average 30 \% (27 \%) of the production, for planner D$^\star$ (S), is stored in the BESS over the entire dataset. $\maxcharge_\%$ is equal to 17.9 \% (17.9 \%\footnote{The value is the same for the $|\Omega|$ and $\sigma$ values considered.}) for the planner D$^\star$ (S) meaning the BESS reached its maximum storage level 5 days out of the 28 days of the dataset. In fact, during sunny days, the BESS is fully charged. A larger BESS capacity should decrease the curtailment and improve the gross revenue. It should be noted that this is a winter month where the maximum generated PV power reached only  half of the installed PV capacity. During a summer month, the maximum production should reach at least 80 \% of the total installed capacity on sunny days. Thus, with a storage capacity of 1 MWh, the curtailment is expected to be much more important during summer sunny days.

\begin{figure}[tb]
	\centering
	\begin{subfigure}{.25\textwidth}
		\centering
		\includegraphics[width=\linewidth]{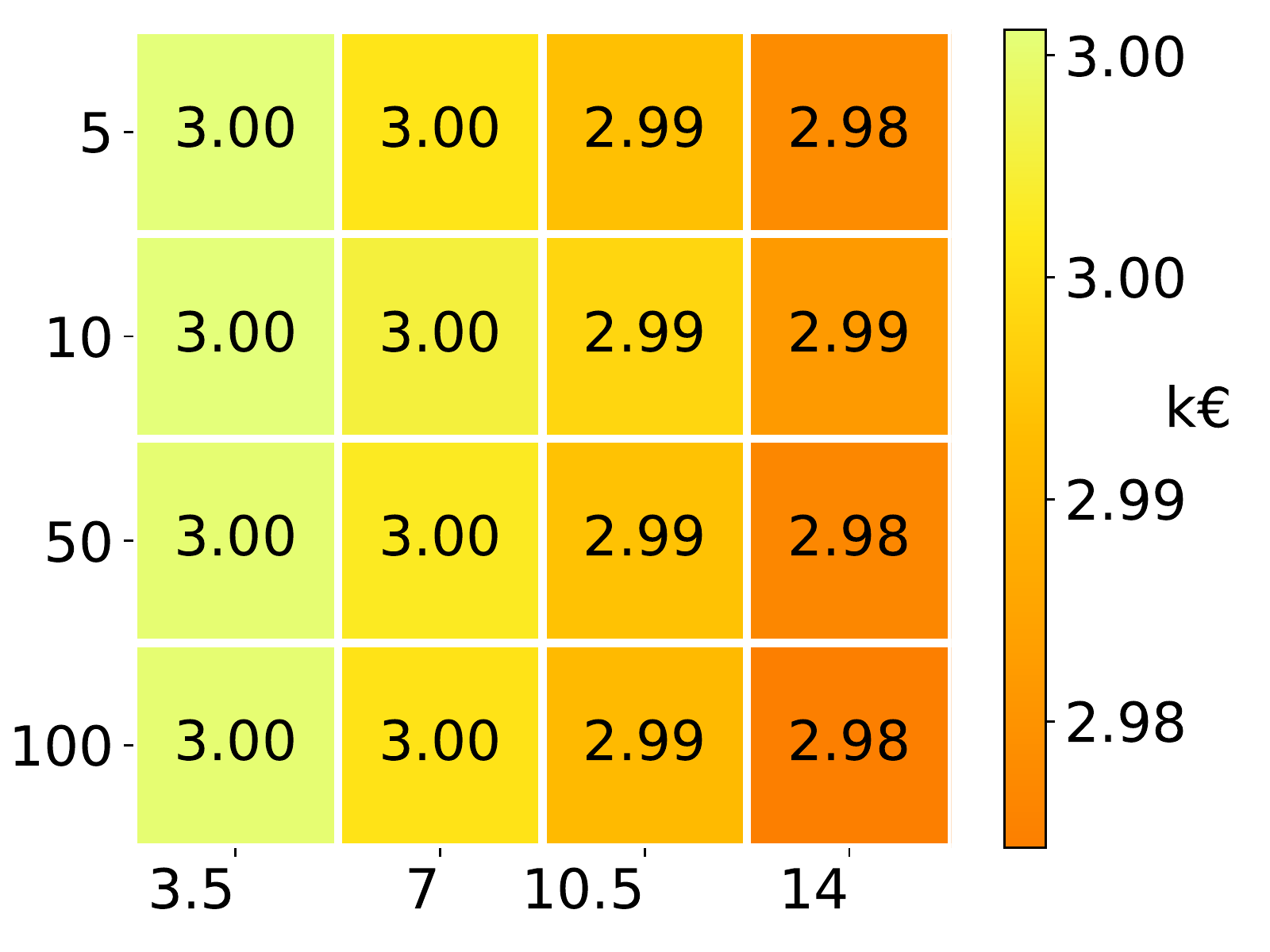}
		\caption{$R^e$.}
	\end{subfigure}%
	\begin{subfigure}{.25\textwidth}
		\centering
		\includegraphics[width=\linewidth]{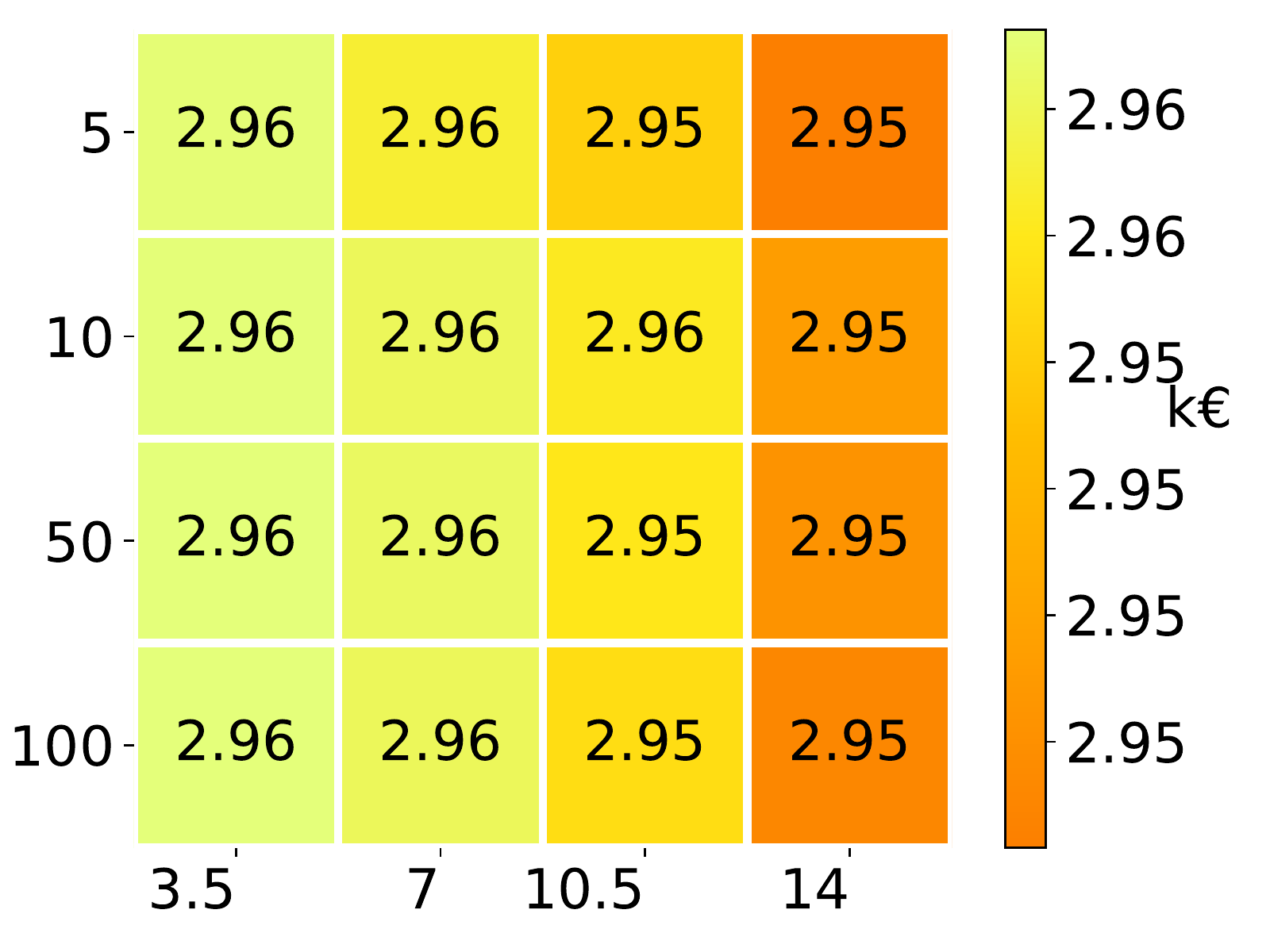}
		\caption{$R^{n,e}=-J^\text{eval}_S$.}
	\end{subfigure}
	\caption{Planner S revenue indicators with $|\Omega|=5, 10, 50, 100$,	$\sigma = 3.5, 7, 10.5, 14\%$.}
	\label{fig:revenue_indicators_planner_S_battery_1000}
\end{figure}

Table~\ref{tab:revenue_indicators_planner_perfect} and Figure~\ref{fig:revenue_indicators_planner_S_battery_1000} provide the results of the revenue indicators for the planners D$^\star$ and S, respectively. It should be noted that in this case, $R^{n,e}=-J^\text{eval}_S$.
\begin{table}[tb]
	\caption{Planner D$^\star$ revenue indicators.}
	\label{tab:revenue_indicators_planner_perfect}
	\centering
	\begin{tabular}{rrrrr} \hline \hline
		$R^e$& $r^e$ & $C^e$ & $R^{n,e}$ & $J^\text{eval}_{D^\star}$	\\ \hline
		3.0 & 94.9 & 0.04  & 2.96 & -2.96  \\ \hline \hline
	\end{tabular}
\end{table}
The smallest value of the objective function is achieved by the planner D$^\star$ and is followed closely by the planner S, even for the highest value of $\sigma$. This result demonstrates the validity of the approach when exploiting an unbiased stochastic predictor. 

In terms of net revenue, both planners achieved 93.7 \% of $R_{\max} \approx 3.16$ k\euro, that results in a loss of 6.3 \%. Most of this loss is due to the curtailment of PV generation. For both planners, the net revenue increases with the generation.

For sunny days, the difference between the nominations and the exports is higher than the deadband just before the production occurs, between 5 and 8 am, and smaller during the main hours of production, between 10 am and 3 pm. Indeed, the planner tends to maximize the revenue by maximizing the exports. However, the ramping power constraints (\ref{eq:D_energy_deviation_constraints}) impose a maximum difference between two consecutive nominations. To maximize the net revenue over the entire day, the planner computes nominations that are not achievable at the beginning of the day to maximize the exports during the day. This results in a penalty, between 5 and 8 am. 

\subsection{BESS capacity sensitivity analysis}

\begin{figure}[!htb]
	\centering
	\begin{subfigure}{.25\textwidth}
		\centering
		\includegraphics[width=\linewidth]{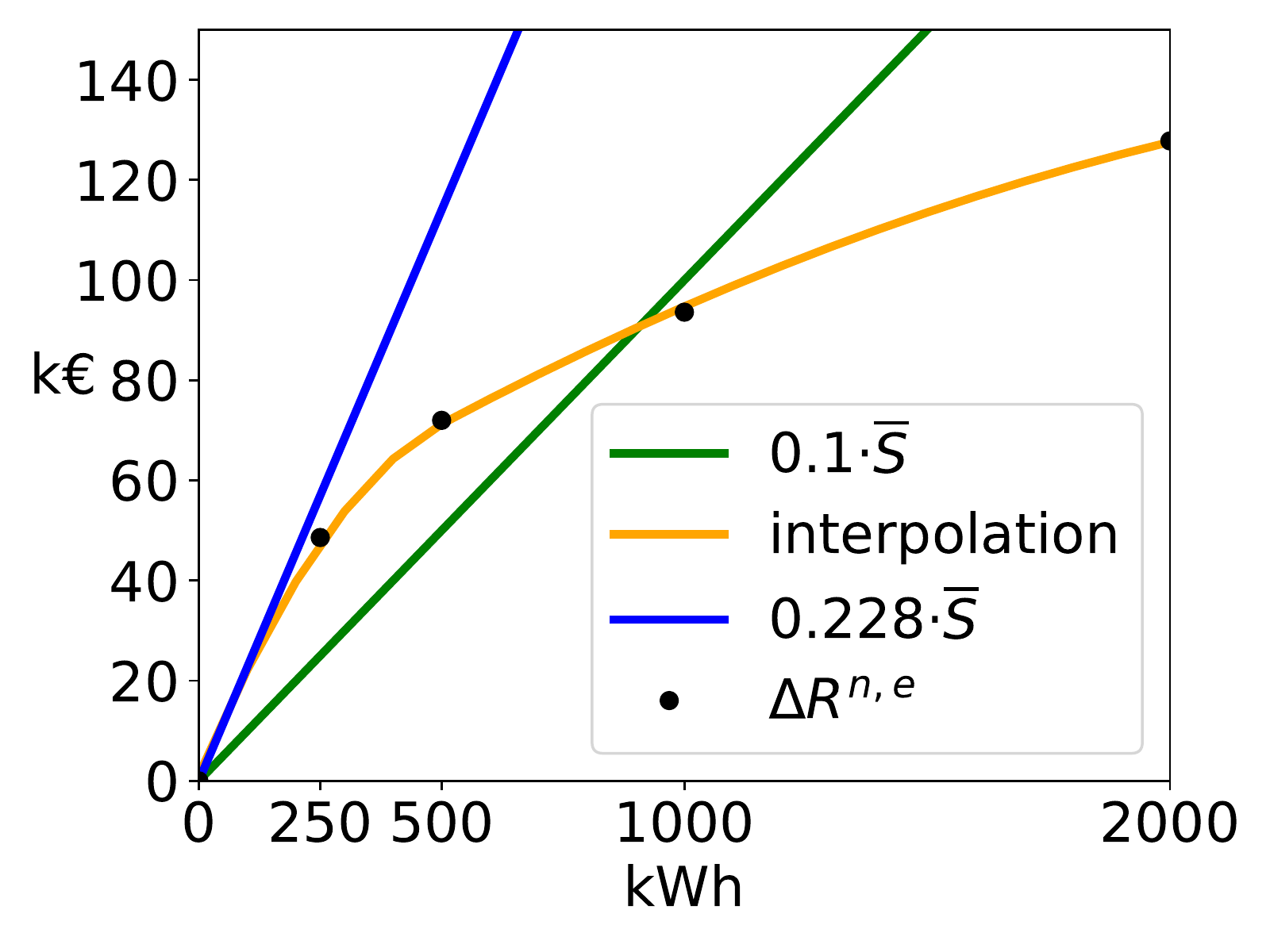}
		\caption{Variation of the net revenue.}
		\label{fig:delta_marginal_gain}
	\end{subfigure}%
	\begin{subfigure}{.25\textwidth}
		\centering
		\includegraphics[width=\linewidth]{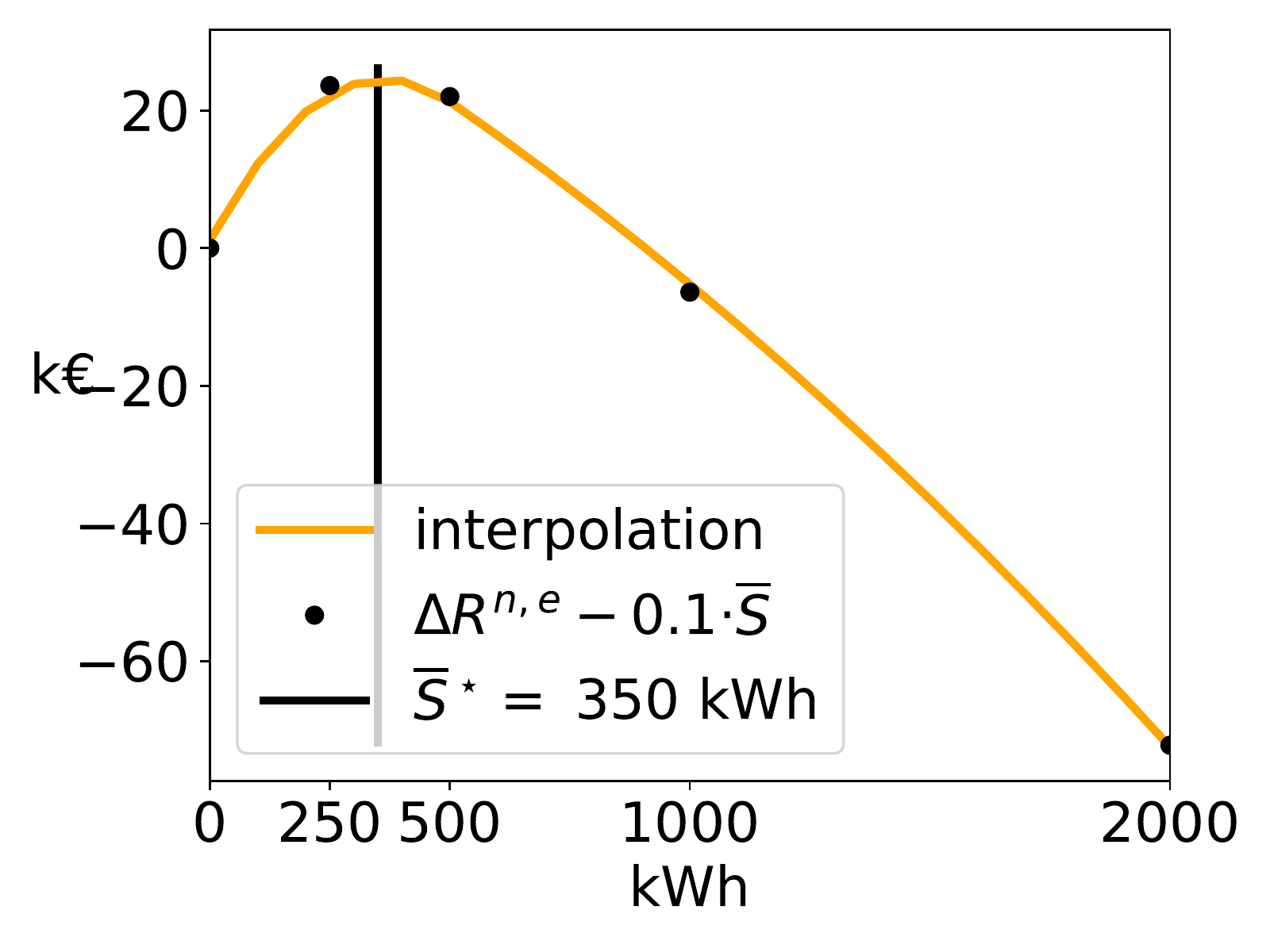}
		\caption{Optimal BESS.}
		\label{fig:optimal_battery_size}
	\end{subfigure}
	\caption{Optimal BESS for a given CAPEX price.}
\end{figure}

The goal is to conduct a sensitivity analysis on the BESS capacity $\maxcharge$ to determine its marginal value and the optimal BESS size $\maxcharge^\star$ for a given CAPEX $\pi_{\maxcharge}$. The efficiencies are still assumed to be unitary. $\initialCharge$ and $\finalCharge$ are set to 0 kWh. Table~\ref{tab:battery_parameters_sensitivity_analysis} provides the other BESS parameters for the five cases.
The scenarios are generated using $|\Omega|=100$, and $\sigma = 3.5, 7, 10.5, 14\%$. 
\begin{table}[tb]
	\caption{BESS parameters.}\label{tab:battery_parameters_sensitivity_analysis}
	\centering
	\begin{tabular}{crrrr} \hline \hline
		Case & $\maxcharge$ [kWh] 	& $\mincharge$ [kWh] & $\chargerate$ [kW]&  $\dischargerate$ [kW]\\ \hline
		1	& 2000 & 0 & 2000 & 2000  \\
		2	& 1000 & 0 & 1000 & 1000  \\
		3	& 500 & 0 & 500 & 500 \\
		4	& 250 & 0 & 250 & 250 \\
		5	& 0 & 0 & 0 & 0 \\ \hline \hline
	\end{tabular}
\end{table}
A new indicator, expressed in k\euro, is defined to quantify the gain provided by the BESS over fifteen years
\begin{align}\label{eq:price_indicator}
& \Delta R^{n,e}_i = 15 \times 12\times (R^{n,e}_i-R^{n,e}_5 ) \quad \forall i \in \{1, 2, 3, 4\}.
\end{align}
It is a lower bound of the total gain as it relies on the results of a winter month. A summer month should provide higher revenue. Table~\ref{tab:battery_sensitivity_capacity_ratio_indicators} provide the planner D$^\star$ indicators. The results demonstrate the interest of using a BESS to optimize the bidding. The larger the BESS is, the lower the curtailment is. Thus, the net revenue increases with the BESS capacity. The maximum achievable revenue is reached with a storage capacity of 2 MWh. However, the larger the BESS is, the smaller $ \Delta R^{n,e}$ increases. It means the marginal benefit decreases with the increase of BESS capacity. A trade-off should be found between the BESS capacity and its CAPEX. Figure~\ref{fig:delta_marginal_gain} provides $ \Delta R^{n,e}$ and its quadratic interpolation in comparison with two BESS prices $\pi_{\maxcharge} = $ 0.1 and 0.228 k\euro/kWh. The value of the derivative $ [\frac{d \Delta R^{n,e}}{d \maxcharge}]_{\maxcharge=0}$ provides the maximum CAPEX that provides a profitable BESS. Then, the optimal storage capacity $\maxcharge^\star$ for a given CAPEX is provided solving $ [\frac{d \Delta R^{n,e}}{d \maxcharge}]_{\maxcharge} = \pi_{\maxcharge}$. For instance, with a CAPEX of 0.1 k\euro/kWh, $\maxcharge^\star$ is approximately 350 kWh. Figure~\ref{fig:optimal_battery_size} provides the values of $  \Delta R^{n,e} - \pi_{\maxcharge} \maxcharge$ with a quadratic interpolation.

\begin{figure}[!htb]
	\centering
	\begin{subfigure}{.25\textwidth}
		\centering
		\includegraphics[width=\linewidth]{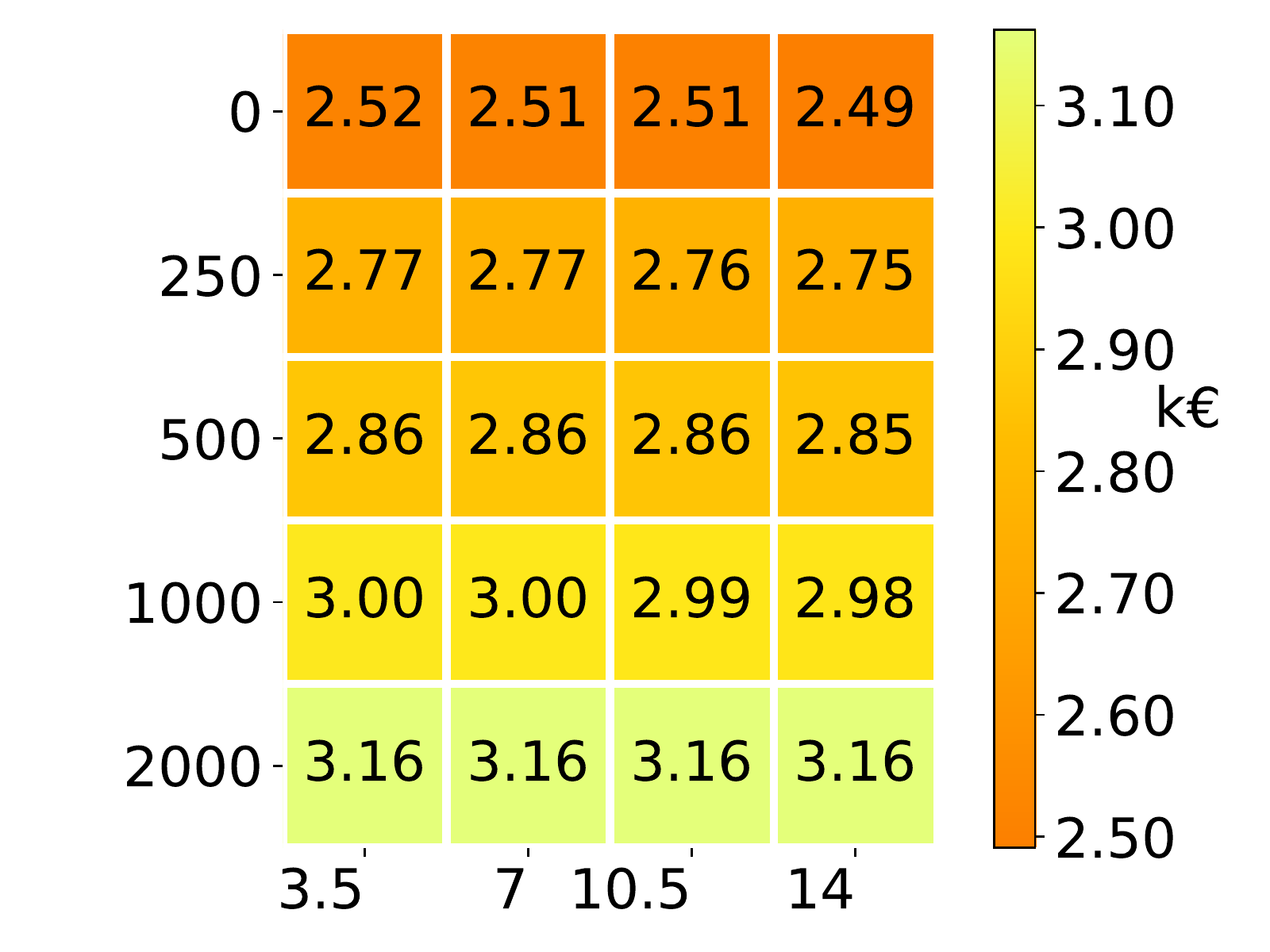}
		\caption{$R^e$.}
	\end{subfigure}%
	\begin{subfigure}{.25\textwidth}
		\centering
		\includegraphics[width=\linewidth]{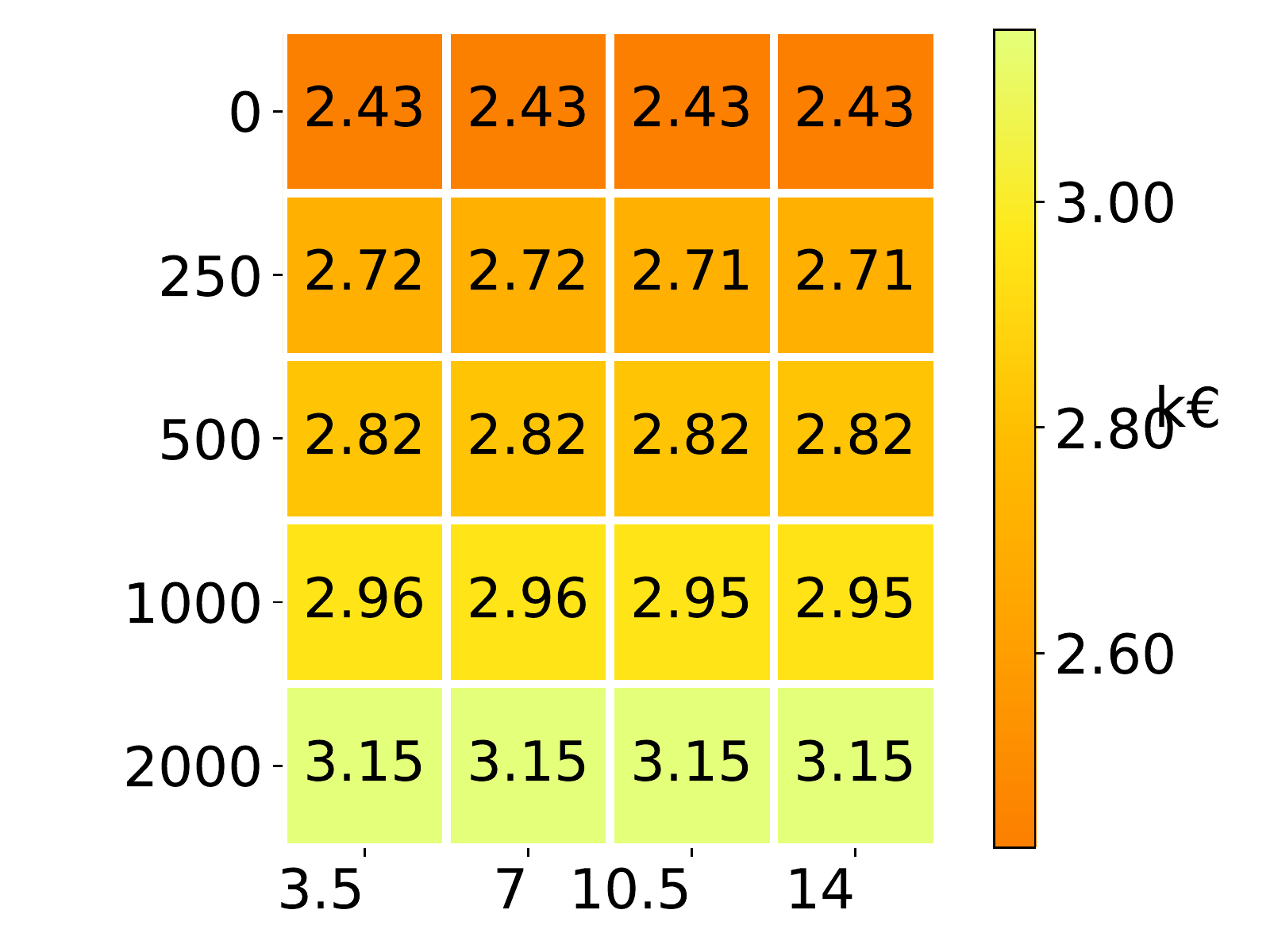}
		\caption{$R^{n,e}=-J^\text{eval}_S$.}
	\end{subfigure}
	\caption{Planner S revenue indicators BESS capacity sensitivity analysis with $|\Omega|=100$, $\sigma = 3.5, 7, 10.5, 14\%$.}
	\label{fig:revenue_indicators_planner_S_battery_battery_sensitivity}
\end{figure}

Figure~\ref{fig:revenue_indicators_planner_S_battery_battery_sensitivity} provides the planner S revenue indicators. The results are still almost identical for all indicators for the smallest value of $\sigma$ and very close with the highest one, as is expected.   
\begin{table}[tb]
	\caption{Planner D$^\star$ ratio and revenue indicators BESS capacity sensitivity analysis.}
	\label{tab:battery_sensitivity_capacity_ratio_indicators}
	\centering
	\begin{tabular}{crrrrr} \hline \hline
		Case   & $[\production]^D$ & 	$\production_\% $& $\charge_\%  $ & $\maxcharge_\%$ & $\export_\% $	\\ \hline
		1 &   70.3 & 	100  & 45.6 & 17.9 & 77.7 \\ 
		2 &   66.7 &	94.9 & 29.6 & 17.9 & 76.2 \\ 
		3 &   63.6 & 	90.5 & 17.3 & 39.3 & 76.4 \\ 
		4 &   61.7 & 	87.7 & 11.1 & 46.4 & 75.2 \\ 
		5 &   56.0 & 	79.6 & -    & -    & 70.4 \\ \hline
		Case    & $R^e$& $C^e$  & $R^{n,e}$ & $J^\text{eval}_{D^\star}$	& $\Delta R^{n,e}$\\ \hline
		1 &  3.16 & 0.01  & 3.15 & -3.15 & 128\\
		2 &  3.0  & 0.04  & 2.96 & -2.96 & 94\\
		3 &  2.86 & 0.04  & 2.84 & -2.84 & 72\\
		4 &  2.77 & 0.06  & 2.71 & -2.71 & 49\\
		5 &  2.52 & 0.08  & 2.44 & -2.44 & 0 \\\hline \hline
	\end{tabular}
\end{table}
%
%
%
%

\section{Conclusions and perspectives}\label{sec:conclusions}

This paper addresses the energy management of a grid-connected PV plant coupled with a BESS within the capacity firming framework, which is decomposed in two steps: computing the day-ahead nominations, then computing the renominations and the set-points in real-time to minimize the energy and ramp power deviations from nominations. This paper investigates the first step by comparing a stochastic and a deterministic formulation. The main goal is to validate the stochastic approach by using an ideal predictor providing unbiased PV scenarios. 

The results of the stochastic planner are comparable with those of the deterministic planner, even when the prediction error variance is non-negligible. Finally, the results of the BESS capacity sensitivity analysis demonstrate the advantage of using a BESS to optimize the bidding day-ahead strategy. However, a trade-off must be found between the marginal gain provided by the BESS and its investment and operational costs.

Several extensions of this work are under investigation.
The first is to better assess the planner's behavior using a full year of data. Then, the next challenge is to use a more realistic methodology to generate PV generation scenarios. Several scenario generation approaches could be investigated, based on a point forecast model such as the PVUSA model \cite{dows1995pvusa,bianchini2013model,bianchini2020estimation}, combined with Gaussian copula \cite{papaefthymiou2008using,pinson2012evaluating,golestaneh2016generation}. 
Another challenge is to consider the non-convex penalty function specified by the CRE into the objective.
Finally, the last challenge is to investigate the second step of the capacity firming problem, for instance by adapting the approach implemented in \cite{dumas2020coordination}.

\bibliographystyle{IEEEtran}
\bibliography{mybib}

\begin{thebibliography}{10}
\providecommand{\url}[1]{#1}
\csname url@samestyle\endcsname
\providecommand{\newblock}{\relax}
\providecommand{\bibinfo}[2]{#2}
\providecommand{\BIBentrySTDinterwordspacing}{\spaceskip=0pt\relax}
\providecommand{\BIBentryALTinterwordstretchfactor}{4}
\providecommand{\BIBentryALTinterwordspacing}{\spaceskip=\fontdimen2\font plus
\BIBentryALTinterwordstretchfactor\fontdimen3\font minus
  \fontdimen4\font\relax}
\providecommand{\BIBforeignlanguage}[2]{{%
\expandafter\ifx\csname l@#1\endcsname\relax
\typeout{** WARNING: IEEEtran.bst: No hyphenation pattern has been}%
\typeout{** loaded for the language `#1'. Using the pattern for}%
\typeout{** the default language instead.}%
\else
\language=\csname l@#1\endcsname
\fi
#2}}
\providecommand{\BIBdecl}{\relax}
\BIBdecl

\bibitem{pinson2007trading}
P.~Pinson, C.~Chevallier, and G.~N. Kariniotakis, ``Trading wind generation
  from short-term probabilistic forecasts of wind power,'' \emph{IEEE
  Transactions on Power Systems}, vol.~22, no.~3, pp. 1148--1156, 2007.

\bibitem{bitar2012bringing}
E.~Y. Bitar, R.~Rajagopal, P.~P. Khargonekar, K.~Poolla, and P.~Varaiya,
  ``Bringing wind energy to market,'' \emph{IEEE Transactions on Power
  Systems}, vol.~27, no.~3, pp. 1225--1235, 2012.

\bibitem{giannitrapani2014bidding}
A.~Giannitrapani, S.~Paoletti, A.~Vicino, and D.~Zarrilli, ``Bidding strategies
  for renewable energy generation with non stationary statistics,'' \emph{IFAC
  Proceedings Volumes}, vol.~47, no.~3, pp. 10\,784--10\,789, 2014.

\bibitem{giannitrapani2015bidding}
A.~Giannitrapani, S.~Paoletti, A.~Vicino, and Zarrilli, ``Bidding wind energy
  exploiting wind speed forecasts,'' \emph{IEEE Transactions on Power Systems},
  vol.~31, no.~4, pp. 2647--2656, 2015.

\bibitem{riffonneau2011optimal}
Y.~Riffonneau, S.~Bacha, F.~Barruel, and S.~Ploix, ``Optimal power flow
  management for grid connected pv systems with batteries,'' \emph{IEEE
  Transactions on sustainable energy}, vol.~2, no.~3, pp. 309--320, 2011.

\bibitem{haessig2014dimensionnement}
P.~Haessig, ``Dimensionnement et gestion d’un stockage d’{\'e}nergie pour
  l'att{\'e}nuation des incertitudes de production {\'e}olienne,'' Ph.D.
  dissertation, Cachan, Ecole normale sup{\'e}rieure, 2014.

\bibitem{haessig2013aging}
P.~Haessig, B.~Multon, H.~B. Ahmed, S.~Lascaud, and L.~Jamy, ``Aging-aware nas
  battery model in a stochastic wind-storage simulation framework,'' in
  \emph{2013 IEEE Grenoble Conference}.\hskip 1em plus 0.5em minus 0.4em\relax
  IEEE, 2013, pp. 1--6.

\bibitem{haessig2015energy}
P.~Haessig, B.~Multon, H.~B. Ahmed, S.~Lascaud, and P.~Bondon, ``Energy storage
  sizing for wind power: impact of the autocorrelation of day-ahead forecast
  errors,'' \emph{Wind Energy}, vol.~18, no.~1, pp. 43--57, 2015.

\bibitem{n2019optimal}
A.~N'Goran, B.~Daugrois, M.~Lotteau, and S.~Demassey, ``Optimal engagement and
  operation of a grid-connected pv/battery system,'' in \emph{2019 IEEE PES
  Innovative Smart Grid Technologies Europe (ISGT-Europe)}.\hskip 1em plus
  0.5em minus 0.4em\relax IEEE, 2019, pp. 1--5.

\bibitem{dows1995pvusa}
R.~Dows and E.~Gough, ``Pvusa procurement, acceptance, and rating practices for
  photovoltaic power plants,'' Pacific Gas and Electric Co., San Ramon, CA
  (United States). Dept. of~…, Tech. Rep., 1995.

\bibitem{bianchini2013model}
G.~Bianchini, S.~Paoletti, A.~Vicino, F.~Corti, and F.~Nebiacolombo, ``Model
  estimation of photovoltaic power generation using partial information,'' in
  \emph{IEEE PES ISGT Europe 2013}.\hskip 1em plus 0.5em minus 0.4em\relax
  IEEE, 2013, pp. 1--5.

\bibitem{bianchini2020estimation}
G.~Bianchini, D.~Pepe, and A.~Vicino, ``Estimation of photovoltaic generation
  forecasting models using limited information,'' \emph{Automatica}, vol. 113,
  p. 108688, 2020.

\bibitem{papaefthymiou2008using}
G.~Papaefthymiou and D.~Kurowicka, ``Using copulas for modeling stochastic
  dependence in power system uncertainty analysis,'' \emph{IEEE Transactions on
  Power Systems}, vol.~24, no.~1, pp. 40--49, 2008.

\bibitem{pinson2012evaluating}
P.~Pinson and R.~Girard, ``Evaluating the quality of scenarios of short-term
  wind power generation,'' \emph{Applied Energy}, vol.~96, pp. 12--20, 2012.

\bibitem{golestaneh2016generation}
F.~Golestaneh, H.~B. Gooi, and P.~Pinson, ``Generation and evaluation of
  space--time trajectories of photovoltaic power,'' \emph{Applied Energy}, vol.
  176, pp. 80--91, 2016.

\bibitem{dumas2020coordination}
J.~Dumas, S.~Dakir, C.~Liu, and B.~Corn{\'e}lusse, ``Coordination of
  operational planning and real-time optimization in microgrids,'' in \emph{XXI
  Power Systems Computation Conference}, 2020.

\bibitem{box2015time}
G.~E. Box, G.~M. Jenkins, G.~C. Reinsel, and G.~M. Ljung, \emph{Time series
  analysis: forecasting and control}.\hskip 1em plus 0.5em minus 0.4em\relax
  John Wiley \& Sons, 2015.

\end{thebibliography}


\section{Annex: PV scenario generation methodology}\label{sec:scenario_generation_methodology}

This Annex describes the methodology to generate the set of unbiased PV scenarios. The goal is to define an ideal unbiased predictor with a fixed variance over all lead times. In this section, let $t$ be the current time index, $k$ be the lead time of the prediction, $K$ be the maximum lead time of the prediction, $y_{t+k}$ be the true value of the signal $y$ at time $t + k$, and $\widehat{y}_{t+k|t}$ be the value of $y_{t+k}$ predicted at time $t$. The forecasts are computed at 4 pm (nominations deadline) for the day-ahead. With a market period duration of fifteen minutes, $K$ is equal to 128. The PV forecasts are needed for lead times from $k=33$ (00:00 to 00:15 am) to $k=K=128$ (11:45 to 12:00 pm). Then, $\widehat{y}_{t+k|t}$ and $y_{t+k}$ are assumed to be related by
\begin{align}\label{eq:prediction_error_definition}	
\widehat{y}_{t+k|t}  & =  y_{t+k} (1+\epsilon_k).
\end{align}
The error term $\epsilon_k$ is generated by the moving-average model defined in Chapter 3 of \cite{box2015time}
\begin{subequations}
	\begin{align}\label{eq:epsilon_definition}	
	\epsilon_1 & =  \eta_1 \\ 
	\epsilon_k & =  \eta_k + \sum_{i=1}^{k-1} \alpha_i \eta_{k-i} \quad \forall k \in \{2, ..., K\}
	\end{align}
\end{subequations}
with $\{ \alpha_i \}_{i=1}^{K-1} $ scalar coefficients, $\{\eta_k \}_{k=1}^K $ independent and identically distributed sequences of random variables from a normal distribution $\mathcal{N}(0, \sigma)$. Thus, the variance of the error term is
\begin{subequations}
	\begin{align}\label{eq:epsilon_var_definition}	
	& \mathrm{Var} (\epsilon_1)  =  \sigma^2 \\ 
	& \mathrm{Var} (\epsilon_k)  =  \big( 1 + \sum_{i=1}^{k-1} \alpha_i^2 \big) \sigma^2 \quad \forall k \in \{2, ..., K\}.
	\end{align}
\end{subequations}
It is possible to simulate with this model an increase of the prediction error variance with the lead time $k$ by choosing
\begin{align}\label{eq:alpha_definition}	
\alpha_i & =  p^i \quad \forall i \in \{1, ..., K-1\}.
\end{align}
(\ref{eq:epsilon_var_definition}) becomes, $\forall k \in \{1, ..., K\}$
\begin{align}\label{eq:epsilon_definition_2}	
\mathrm{Var} (\epsilon_k) &  = \sigma^2 A_{\epsilon_k} 
\end{align}
with $A_{\epsilon_k}$ defined $\forall k \in \{1, ..., K\}$ by
\begin{align}\label{eq:Ak_definition}	
A_{\epsilon_k} &  =  \sum_{i=0}^{k-1} (p^2)^i = \frac{1-(p^2)^{k}}{1-p^2}.
\end{align}
Then, with $0 \leq p < 1 $, it is possible to make the prediction error variance independent of the lead time as it increases. Indeed:
\begin{align}\label{eq:Ak_limit}	
\lim_{k \to\infty}  A_{\epsilon_k} & = A_\infty = \frac{1}{1-p^2}.
\end{align}
For instance, with $p=0.9$ and $K=128$, for $k \geq 33$, $A_{\epsilon_k}\approx A_\infty$ that is approximately 5.26. Thus, $\forall k \in \{33, ..., K\}$
\begin{align}\label{eq:epsilon_variance_value}	
& \mathrm{Var} (\epsilon_k)  \approx 5.26  \sigma^2.
\end{align}
Finally, the $\sigma$ value to set a maximum $\epsilon_\text{max}$ with a high probability of 0.997, corresponding to a three standard deviation confidence interval from a normal distribution, is found by imposing $\epsilon_\text{max} = 3 \sqrt{\mathrm{Var} (\epsilon_K)}$:
\begin{align}\label{eq:epsilon_max_value}	
\sigma & \approx \frac{\epsilon_\text{max}}{3\sqrt{A_\infty}},
\end{align}
with $\epsilon_\text{max} = 0.25, 0.50, 0.75, 1.0$, $\sigma = 3.5, 7, 10.5, 14\%$.

\end{document}